\documentclass{amsproc}

\newtheorem{theorem}{Theorem}[section]
\newtheorem{lemma}[theorem]{Lemma}

\theoremstyle{definition}
\newtheorem{definition}[theorem]{Definition}
\newtheorem{example}[theorem]{Example}
\newtheorem{xca}[theorem]{Exercise}

\theoremstyle{remark}
\newtheorem{remark}[theorem]{Remark}

\numberwithin{equation}{section}

\theoremstyle{corollary}
\newtheorem{corollary}[theorem]{Corollary}

\newcommand{\abs}[1]{\lvert#1\rvert}

\newcommand{\blankbox}[2]{%
  \parbox{\columnwidth}{\centering
    \setlength{\fboxsep}{0pt}%
    \fbox{\raisebox{0pt}[#2]{\hspace{#1}}}%
  }%
}

\newcommand{\Z}{{\mathbf Z}}

\newcommand{\R}{{\mathbf R}}
\newcommand{\C}{{\mathbf C}}

\newcommand{\id}{{\rm {id }}}

\newcommand{\tc}{{\rm\bf {TC}}}

\newcommand{\comment}[1]{}

\begin{document}

\title[Moving obstacles]{Topological complexity of collision free motion planning
algorithms in the presence of multiple moving obstacles}

\author{Michael Farber}
\address{Department of Mathematics, University of Durham, Durham DH1 3LE, UK}
\email{Michael.Farber@durham.ac.uk}

\author{Mark Grant}
\address{Department of Mathematics, University of Durham, Durham DH1 3LE, UK}

\email{Mark.Grant@durham.ac.uk}


\author{Sergey Yuzvinsky}
\address{Department of Mathematics, University of Oregon, Eugene}
\email{yuz@math.uoregon.edu}
\date{January 1, 1994 and, in revised form, June 22, 1994.}



\begin{abstract}
We study motion planning algorithms for collision free control of
multiple objects in the presence of moving obstacles. We compute
the topological complexity of algorithms solving this problem. We
apply topological tools and use information about cohomology
algebras of configuration spaces. The results of the paper may
potentially be used in systems of automatic traffic control.
\end{abstract}

\maketitle
\section{Introduction}

The theory of robot motion planning \cite{L}, \cite{HS} has developed a broad
variety of algorithms designed for different real life situations. Most of
these algorithms assume that the robot performs an online search of the scene
and eventually finds its way to the goal. In this approach one understands that
the robot has initially a very limited knowledge of the scene and supplements
this knowledge by using sensors, vision, memory and perhaps some ability to
analyze these data.

A completely different situation arises when one controls
simultaneously multiple objects (robots) moving in a coordinated
way in a fully known environment.  Assume, for instance, that we
have $n$ objects moving in $\R^3$ with no collisions and avoiding
the obstacles whose geometry is prescribed in advance. In this
case the dimension $3n$ of the configuration space of the system
is large (if the number $n$ of controlled objects is large) and
therefore the online search algorithms as described above become
much less effective.

The topological approach to the robot motion planning problem initiated in
\cite{F1, F2} is applicable in some situations with high dimensional
configuration space under the assumption that the configuration space of the
system is known in advance. One divides the whole configuration space into
pieces (local domains) and prescribes continuous motions (local rules) over
each of the local domains. The minimal number of such local domains is the
measure of topological complexity of the problem.

Formally, the topological complexity is a numerical invariant
$\tc(X)$ of the configuration space of the system \cite{F1}. Its
introduction was inspired by the earlier well-known work of S.
Smale \cite{Sm} and V. Vassiliev \cite{Va} on the theory of
topological complexity of algorithms of solving polynomial
equations. The approach of \cite{F1, F2} was also based on the
general theory of robot motion planning algorithms described in
the book of J.-C. Latombe \cite{L} and on the abstract theory of a
genus of a fiber space developed by A.S. Schwarz \cite{Sz}. Paper
\cite{F4} contains a recent survey.

Paper \cite{FY} solves the problem of finding the topological
complexity of motion planning algorithms for controlling many
particles moving in $\R^3$ or on the plane $\R^2$ with no
collisions. Such motion planning algorithms appear in automatic
traffic control problems. Paper \cite{FY} uses the techniques of
the theory of subspace arrangements \cite{OT}. It was shown in
\cite{FY} that the complexity is roughly twice the number of
particles (it is $2n-2$ for the planar case and $2n-1$ for the
spatial case). At the moment we do not know specific motion
planning algorithms with complexity linear in $n$. However a
quadratic in $n$ motion planning algorithm was described in \S 26
of \cite{F4}.

In paper \cite{F3} we studied the problem of computing the
topological complexity of the motion planning problem for $n$
particles moving on a graph with forbidden collisions. It was
shown that for large numbers of particles the topological
complexity depends only on the geometry of the graph and is
independent of the number of particles. Some specific motion
planning algorithms for this problem were also constructed
\cite{F3}.

In the present paper we study algorithms solving the following motion planning
problem. Several objects are to be transported from an initial configuration
$A_1, A_2, \dots, A_n\in \R^3$ to a final configuration $B_1, B_2, \dots,
B_n\in \R^3$ such that in the process of motion there occur no collision
between the objects and such that the objects do not touch the obstacles in the
process of motion.

We make several simplifying assumptions: (a) Each object is
represented by a single point; (b) each obstacle is represented by
a single point; (c) collision between two objects occurs if they
are situated at the same point in space; (d) an object touches an
obstacle if the points representing the object and the obstacle
coincide. We also make an important assumption that the behavior
of the obstacles is known in advance. We show that because of this
assumption the problem becomes topologically equivalent to the
similar problem when the obstacles are stationary.

One of our results seems surprising: we show that complexity of collision free
navigation of many objects in the presence of moving obstacles is essentially
independent of the number of obstacles and grows linearly with the number of
controlled objects.

We believe that the conclusions of this paper will remain valid in a more
general and realistic situation when the objects and the obstacles are
represented by small balls, possibly of different radii, and the control
requirements are to avoid tangencies between objects and obstacles.

\section{Statement of the result}

Let us explain our notation. $n$ denotes the number of moving objects and $m$
the number of obstacles. The symbol $F(X, n)$ denotes the configuration space
of $n$ distinct points in a topological space $X$; in other words, $F(X,n)$ is
the subset of the Cartesian power $X\times X\times \dots \times X$ ($n$ times)
consisting of configurations $(x_1, \dots, x_n)$ (where $x_i\in X$) with
$x_i\not=x_j$ for $i\not=j$.

We assume that the controllable objects lie in $\R^{r+1}$ where $r\geq 1$ is an
integer. Clearly, the practical applications are covered by the cases $r=1$ or
$r=2$ which will be considered in more detail later. In this section we will
assume that $r$ is arbitrary.

The algorithms we study take as input the following data:
\begin{enumerate} \item[(a)] the initial configuration $A=(A_1,
A_2, \dots, A_n)\in F(\R^{r+1}, n)$ of the objects, $A_i\in
\R^{r+1}$, with $A_i\not= A_j$; \item[(b)] the desired final
configuration of the objects $B=(B_1, B_2, \dots, B_n)\in
F(\R^{r+1}, n)$, where $B_i\in \R^{r+1}$, $B_i\not= B_j$;
\item[(c)] the trajectory of moving obstacles $C(t)=(C_1(t),
C_2(t), \dots, C_m(t))$ $\in F(\R^{r+1}, m)$, $C_i(t)\not=C_j(t)$
for any $t\in [0,1]$.
\end{enumerate}
Here $t\in [0,1]$ is the time. We assume that $A_i\not= C_j(0)$ and
$B_i\not=C_j(1)$ for all $i=1, \dots, n$ and $j=1, \dots, m$. This means that
initially the objects are all distinct and disjoint from the obstacles at time
$t=0$ and the desired positions of the objects are pairwise distinct and
distinct from positions of the obstacles at time $t=1$.

The output of the algorithm is a continuous motion of the objects
$$\gamma_i(t)\in \R^{r+1}, \quad t\in [0,1], \quad i=1, \dots, n$$ satisfying the
following conditions: \begin{enumerate} \item[($\alpha$)] $\gamma_i(t)$ are
continuous functions of $t$; \item[($\beta$)] $\gamma_i(0)=A_i$ and
$\gamma_i(1)=B_i$; \item[($\gamma$)] for any $t\in [0,1]$ one has
$\gamma_{i}(t)\not=\gamma_{k}(t)$ for $i\not=k$ where $i, k\in \{1, \dots,
n\}$; \item[($\delta$)] for any $t\in [0,1]$ one has $\gamma_i(t) \not= C_j(t)$
where $i=1, \dots, n$ and $j=1, \dots, m$.
\end{enumerate}

Let the function $$Y_m: [0,1]\to F(\R^{r+1}, m)$$ represent the motion of the
obstacles. For any time $t\in [0,1]$ the obstacles are positioned at the points
$$Y_m(t)=(C_1(t), C_2(t), \dots, C_m(t)).$$
Next we introduce a few more notations:

$$E(Y_m) = \{\gamma: [0,1]\to F(\R^{r+1}, n);\,\,
\gamma(t)\cap Y_m(t)=\emptyset, \quad\forall t\in [0,1]\},$$

$$B(Y_m) = \{(A,B)\in F(\R^{r+1},n)\times F(\R^{r+1}, n); A\cap Y_m(0)=\emptyset, \quad
B\cap Y_m(1)=\emptyset\}.$$ Elements of $E(Y_m)$ are motions of the objects
avoiding the obstacles, elements of $B(Y_m)$ represent the input data: the
initial configuration $A$ and the desired configuration $B$.

There is a canonical map
\begin{eqnarray}\label{proj}
\pi(Y_m): E(Y_m) \to B(Y_m)\end{eqnarray} given by
$\pi(Y_m)(\gamma) = (\gamma(0), \gamma(1))$. A {\it section} of
$\pi$ is a map (possibly discontinuous)
\begin{eqnarray}\label{section}
s: B(Y_m)\to E(Y_m),
\end{eqnarray}
satisfying $\pi(Y_m)\circ s=1_{B(Y_m)}.$ Clearly, $s$ assigns to any input data
$(A,B)\in B(Y_m)$ a motion of the controllable objects starting at
configuration $A$, ending at $B$ and avoiding mutual collisions between objects
and between objects and the moving obstacles. We conclude that algorithms
solving our motion planning problem are in one-to-one correspondence with
sections of map (\ref{proj}) (we will see below that it is a fibration). Given
a section $s$ as above, then for any pair $(A,B)\in B(Y_m)$, the output of the
algorithm is $s(A,B)\in E(Y_m)$.

Hence, our task is to estimate the complexity of finding a section of
$\pi(Y_m)$.

We recall that {\it the Schwarz genus} \cite{Sz} of a continuous map $f: X\to
Y$ map is the minimal integer $k$ such that $Y$ can be covered by $k$ open sets
$U_1\cup U_2\dots \cup U_k=Y$ with the property that over each $U_i$ there
exists a continuous section $s_i: U_i\to X$, $f\circ s_i = 1_{U_i}$. We apply
this notion to estimate the complexity of finding a section to (\ref{proj}).

\begin{definition}\label{def1}
The complexity of the motion planning problem for moving $n$
objects in the presence of $m$ moving obstacles is defined as the
Schwarz genus of map (\ref{proj}).
\end{definition}

In practical terms this means the following. One divides the space of all
possible inputs $B(Y_m)$ into several domains (called {\it the local domains})
and specifies a continuous motion planning algorithm (called {\it a local
rule}) over each of the domains. The complexity of the problem is defined
(according to Definition \ref{def1}) as the minimal number of local domains in
all possible algorithms of this type.

Now we state our main result:

\begin{theorem}\label{main}
The complexity of motion planning problem of moving $n\geq 2$ objects in $\R^3$
avoiding collisions in the presence of $m\geq 1$ moving obstacles equals
$2n+1$. The complexity of the similar planar problem {\rm (}i.e. when the
objects are restricted to lie in the plane $\R^2${\rm )} is $2n$ if $m=1$ and
it is $2n+1$ for $m\geq 2$.
\end{theorem}

This theorem follows by combining Theorem \ref{thm31} and Theorems
\ref{thm51} and \ref{thm61} which are stated and proven below.

\section{Reduction to the case of stationary obstacles}

The following theorem shows that Definition \ref{def1} gives the
notion of complexity which coincides with a special case of the
notion of navigational complexity of topological spaces $\tc(X)$
which was studied previously \cite{F1}. Recall that $\tc(X)$
denotes the Schwarz genus of the path space fibration
\begin{eqnarray}\label{tcmap}
\pi: X^I\to X\times X \end{eqnarray} where $X^I$ denotes the space of all
continuous paths $\gamma: I=[0,1]\to X$ with the compact-open topology and
$\pi(\gamma)=(\gamma(0),\gamma(1))$.

\begin{theorem}\label{thm31}
(a) The map $\pi(Y_m)$ (given by (\ref{proj})) is a locally
trivial fibration; (b) The fiberwise homeomorphism type of
fibration $\pi(Y_m)$ is independent of the trajectory of moving
obstacles $Y_m$; in particular, it equals the Schwarz genus of the
special case of (\ref{proj}) when the obstacles are stationary;
(c) The Schwarz genus of $\pi(Y_m)$ equals
$$\tc(F(\R^{r+1}-S_m, n))$$ where $S_m=\{s_1, \dots, s_m\}\subset \R^{r+1}$ is an arbitrary fixed
$m$-element subset.
\end{theorem}
\begin{proof} First we find a continuous family of homeomorphisms $$\psi_t:
\R^{r+1}\to \R^{r+1}, \quad t\in [0,1]$$ having the properties
\begin{eqnarray}\label{equ}
\psi_0=\id, \quad \mbox{and}\quad \psi_t(Y_m(t)) = Y_m(0), \quad t\in [0,1].
\end{eqnarray}
Equation (\ref{equ}) is understood as an equality between ordered
sets; in other words,  (\ref{equ}) states that \begin{eqnarray}
\psi_t(C_i(t))=C_i(0), \quad t\in [0,1].\end{eqnarray} The
existence of such isotopy $\psi_t$ is a standard fact of the
manifold topology; formally it follows from the well known isotopy
extension theorem.

Next we consider the commutative diagram of continuous maps
\begin{eqnarray}\label{diag}
\begin{array}{ccc}
E(Y_m) & \stackrel{F}\to & X^I\\ \\
{\quad \quad \downarrow \pi(Y_m)} & & \downarrow \pi\\ \\
B(Y_m) & \stackrel G\to & X\times X
\end{array}
\end{eqnarray}
where $$X=F(\R^{r+1}-Y_m(0), n)$$ and $\pi$ is the canonical path
space fibration (\ref{tcmap}) for the special case
$X=F(\R^{r+1}-Y_m(0), n)$. Maps $F$ and $G$ are defined as
follows:
\begin{eqnarray}
F(\gamma)(t)= \psi_t(\gamma(t)), \quad \gamma\in E(Y_m), \quad t\in [0,1],\\ \nonumber \\
G(A,B) = (\psi_0(A), \psi_1(B)), \quad (A,B)\in B(Y_m).
\end{eqnarray}
Clearly, $F$ and $G$ are homeomorphisms and diagram (\ref{diag})
is commutative. This implies statements (a) - (c) of Theorem
\ref{thm31}.
\end{proof}

\section{Cohomology of configuration space $F(\R^{r+1}-S_m,n)$}

In this section we collect some known topological results which
will be used in this paper. All cohomology groups have $\Z$ as
coefficients.

 Let
$S_m=\{s_1,\ldots ,s_m\}\subset \R^{r+1}$ be a fixed sequence of
$m$ distinct points in Euclidean space $\R^{r+1}$. Here $r\geq 1$
is an integer. The space $F(\R^{r+1}-S_m,n)$ represents
configurations of $n$ particles in $\R^{r+1}$ avoiding mutual
collisions and collisions with $m$ stationary obstacles
$s_1,\ldots ,s_m$.

The theorems we state in this section are known, they can be found
in Chapter V of the book by Fadell and Husseini \cite{FH}; we have
adjusted the notations as required for our needs. The space
$F(\R^2-S_m, n)$ (here $r=1$) can be identified with the
complement of an arrangement of affine hyperplanes in $\C^n$, so
one can also use the book \cite{OT} as the reference.

We begin by noting that $F(\R^{r+1}-S_m,n)$ embeds in the
configuration space $F(\R^{r+1},n+m)$ via the map $(y_1,\ldots
,y_n)\mapsto (y_1,\ldots ,y_n,s_1,\ldots ,s_m)$ (in fact this
embeds $F(\R^{r+1}-S_m,n)$ as the fibre over $(s_1,\ldots,s_m)$ of
the locally trivial fibration $F(\R^{r+1},n+m)\to F(\R^{r+1},m)$
which projects to the last $m$ points of the configuration). The
cohomology algebra $H^*(F(\R^{r+1},n+m);\Z)$ is well known, and
can be described as follows.

For each pair of integers $i\neq j$ such that $1\leq i,j\leq n+m$,
consider the map $$\phi_{ij}: F(\R^{r+1}, n+m)\to S^{r}, \quad
(y_1, y_2, \ldots, y_{n+m})\mapsto \frac{y_i-y_j}{|y_i-y_j|}\, \in
\, S^{r}.$$  Letting $[S^r]\in H^r(S^r)$ denote the fundamental
cohomology class of the sphere $S^r$, we obtain cohomology classes
$$e_{ij}=\phi_{ij}^\ast[S^r]\in H^r(F(\R^{r+1}, n+m);\Z).$$
\begin{theorem}\label{cohomology}(\cite{FH}, Theorem V.4.2)$\,$
 $H^*(F(\R^{r+1},n+m))$ is the free associative
  graded-commutative algebra generated by the classes $e_{ij}$ for $1\leq i<j\leq n+m$, subject to the relations \\
(i)\quad $e_{ij}^2=0,$ and\\
(ii)\quad $e_{ij}e_{ik}=(e_{ij}-e_{ik})e_{jk}$ for any triple $i<j<k$.
\end{theorem}
\begin{theorem}(\cite{FH}, V.4.2)$\,$ The homomorphism $$H^*(F(\R^{r+1},n+m);\Z)\to
  H^*(F(\R^{r+1}-S_m,n);\Z)$$ induced by
inclusion is an epimorphism, with kernel equal to the ideal
generated by those $e_{ij}$ having $i>n$ and $j>n$.
\end{theorem}
\begin{corollary}\label{relation}
In $H^\ast(F(\R^{r+1}-S_m, n);\Z)$ there are relations
  $$ e_{ij}e_{ik}=0$$
for any triples $i, j, k$ such that $i\leq n$ and $j,k>n$.
\end{corollary}

\begin{theorem}\label{basis}
(\cite{FH}, Theorem V.4.3)$\,$ A basis for
$H^*(F(\R^{r+1}-S_m,n);\Z)$ is given by the set of monomials
$$e_{i_1j_1}e_{i_2j_2}\dots e_{i_sj_s}$$ where $1\leq i_1<i_2<\dots
<i_s\leq n$ and for each $q$ with $1\leq q\leq s$ we have $i_q<j_q$ and $1\leq
j_q\leq n+m$.
\end{theorem}
We remark here that the maximum length of such a monomial is
$(n-1)$ when $m=0$, and $n$ when $m\geq 1$. Hence the highest
dimension in which the integral cohomology of $F(\R^{r+1}-S_m,n)$
is non-trivial is $r(n-1)$ when $m=0$ and $rn$ when $m\geq 1$.
Note also that when $m>1$ the given basis includes all monomials
of the form \begin{eqnarray}\label{bas} \mu_I=\prod_{i\in
I}e_{i,n+1}\prod_{i\not\in I}e_{i,n+2}\end{eqnarray} for every
$I\subset \{1,\ldots,n\}$. For $r=1$ this fact also follows from
the notion of {\bf nbc}-basis for an appropriate ordering of the
generators (see \cite{OT}).

In what follows we will specialise to the cases $r=2$ and $r=1$.

\section{Computing $\tc(F(\R^3-S_m,n))$}
\begin{theorem}\label{thm51} One has:
\begin{displaymath}
\tc(F(\R^3-S_m,n)) = \left\{
\begin{array}{lll}
2n-1& \mbox{if}& m=0,\\
2n+1& \mbox{if}& m\geq 1.
\end{array}\right.
\end{displaymath}
\end{theorem}
\begin{proof}
We first establish a lower bound. We shall use the cohomological
lower bound given by Theorem 7 of \cite{F1}. Set $\bar e_{ij}
=1\otimes e_{ij}-e_{ij}\otimes 1.$ It is a zero-divisor of the
cohomology algebra. Note that $(\bar e_{ij})^2 = -2(e_{ij}\otimes
e_{ij}) \not= 0.$ Consider the following product
$\pi=\prod_{i=1}^{n-1}\, (\bar e_{in})^2.$ We find
$\pi=(-2)^{n-1}\mu\otimes \mu$, where
$\mu=\prod_{i=1}^{n-1}e_{in}.$ The monomial $\mu$ is nonzero by
Theorem \ref{basis}, and hence the product $\pi$ of length
$2(n-1)$ is nonzero. This gives the lower bound
$\tc(F(\R^3-S_m,n))\geq 2n-1$. Now assuming that $m\geq 1$ we have
a nontrivial product $\prod_{i=1}^n\, (\bar e_{i(n+1)})^2$ of
length $2n$, which gives the lower bound $\tc(F(\R^3-S_m,n))\geq
2n+1$ when $m\geq 1$.

To obtain the upper bound, note that $F(\R^3-S_m,n)$ can be viewed
as the complement of a finite collection of codimension 3 affine
subspaces in $\R^{3n}$, so it is simply-connected by an easy
transversality argument. Since it has finitely generated
torsion-free homology and cohomology, it has the homotopy type of
a CW-complex consisting of one $k$-cell for each $k$-dimensional
element in the basis for cohomology given by Theorem \ref{basis}
(see Chapter 4.C of \cite{Ha}). This minimal cell structure is
made explicit in Theorem VI.8.2 of \cite{FH}.  In particular, it
has the homotopy type of a polyhedron of dimension $2(n-1)$ when
$m=0$ and $2n$ when $m\geq 1$. We now apply Corollary 5.3 of
\cite{F2} stating that for a 1-connected polyhedron $Y$,
$$\tc(Y)\leq\mathrm{dim}(Y)+1,$$
which together with homotopy invariance of $\tc$ completes the proof.

\end{proof}

\section{Computing $\tc(F(\R^2-S_m,n))$}

\begin{theorem}\label{thm61} One has:
\begin{displaymath}
\tc(F(\R^2-S_m,n)) = \left\{
\begin{array}{lll}
2n-2& \mbox{if}& m=0,\\
2n& \mbox{if}& m= 1,\\
2n+1& \mbox{if}& m\geq 2.
\end{array}\right.
\end{displaymath}
\end{theorem}

\begin{proof}
The statement for the first two cases follows immediately from \cite{FY} and
the fact that $F(\R^2-X_1,n)$ is homotopy equivalent to $F(\R^2,n+1)$ (see
\cite{FH},p.15).

In what follows we assume that $m\geq 1$. We first establish a
lower bound using again the cohomological lower bound given by
Theorem 7 from \cite{F1}. Set as above $\bar e_{ij} =1\otimes
e_{ij}-e_{ij}\otimes 1$ and consider the product
$$\pi=\prod_{i=1,\ldots,n,\ j=n+1,n+2}\bar e_{ij}.$$
It is clear that $\pi$ can be expressed as a linear combination of
pure tensors $\mu_1\otimes\mu_2$ where $\mu_i$ are monomials in
the $e_{ij}$ complementary to each other in $\pi$. Since the
highest non-zero dimension in $H^*(F(\R^2-S_m,n);\Z)$ is $n$ both
$\mu_i$ should have degree $n$ in order for $\mu_1\otimes\mu_2$
not to vanish. Also it follows from the relation of Corollary
\ref{relation} that the nonvanishing summands of $\pi$ are of the
form $\mu_I\otimes\mu_{\bar I}$, where $\mu_I$ is the monomial
defined by (\ref{bas}),
 $I$ runs over all subsets of
$\{1,2,\ldots,n\}$ and $\bar I$ denotes the complement of $I$.
Since the set $\{\mu_I|I \subset\{1,\ldots,n\}\}$ is a subset of a
basis of $H^n(F(\R^2-S_m,n);\Z)$ given by Theorem \ref{basis}, no
cancellations are possible, whence $\pi\not=0$. This gives the
inequality
 $\tc(F(\R^2-S_m,n))\geq 2n+1$.

The opposite inequality follows immediately from Theorem 5.2 of \cite{F2} since
$F(\R^2-S_m,n)$ has the homotopy type of a connected polyhedron of dimension
$n$ (see \cite{OT}). This completes the proof of the Theorem.
\end{proof}

\section{Concluding remarks}

First we note that our main result Theorem \ref{main} follows by
simply combining Theorems \ref{thm31}, \ref{thm51} and
\ref{thm61}.

Let us compare the following two control problems: (1) motion
planning for moving $n$ objects in $\R^3$ with no collisions and
avoiding collisions with $m\geq 1$ moving obstacles; and (2)
motion planning for moving $n$ objects in $\R^3$ with the only
requirement that they avoid a single point obstacle (in
particular, the objects are allowed to collide, i.e. to occupy the
same position in space). Note that the integer $m$ in problem (1)
can be arbitrarily large. According to Theorem \ref{main} the
topological complexity of the problem (1) is $2n+1$. Surprisingly,
problem (2) also has complexity $2n+1$. Indeed, the configuration
space of problem (2) is Cartesian power of $n$ copies of
$\R^3-\{0\}$ which is homotopy equivalent to
$$S^2\times S^2\times \dots\times S^2\quad (n \, \mbox{times})$$
which has topological complexity $2n+1$ as it is easy to see.
Hence, surprisingly problem (1) which is intuitively more \lq\lq
complicated\rq\rq\, has the same topological complexity as problem
(2).

This comparison shows that in general the notion of topological
complexity is only a partial reflection of real difficulty of a
motion planning problem.


\end{document}